\documentclass[12pt,twoside]{report}
\usepackage{svcon2e}
   \usepackage{amsthm,amssymb}
\usepackage{amsbsy,amsfonts,amsmath}
\usepackage[dvips]{graphics,color}

\catcode`@=12
\makeatletter
\@addtoreset{equation}{section}

\makeatother

\newcommand\Remarks{\medskip\noindent{\bf  Remarks.}\enspace}
\newcommand\Ack{\medskip\noindent{\bf Acknowledgment.}\enspace}

  \newcommand{\Proofmark}[1]{} 

\theoremstyle{plain}

\newtheorem{@definition}{\bf Definition}[section]

\newtheorem{@example}{\bf Example}[section]

\newtheorem{@remark}{\bf Remark}[section]

\newcommand{\nc}{\newcommand}

\nc{\z}{\zeta}
\nc{\ps}{\psi}
\nc{\inv}{^{-1}}
\nc{\ep}{\varepsilon}
\nc{\e}{\epsilon}
\def\ra{\rightarrow}

\nc{\si}{\sigma}
\def\iy{\infty}

\def\be{\begin{equation}}
\def\ee{\end{equation}}
\def\ba{\begin{eqnarray*}}
\def\ea{\end{eqnarray*}}
\def\bae{\begin{eqnarray}}
\def\eae{\end{eqnarray}}
\def\bc{\begin{center}}
\def\ec{\end{center}}
\def\ov{\over}

\def\s{\sigma}

\def\ve{\varepsilon}

\def\p{\Phi_r}
\nc{\noi}{\noindent}

\def\la{\lambda}
\def\pr{\textrm{Prob}}

\def\dl{\delta}

\def\t{\tau}
\def\be{\beta}

\nc{\ga}{\gamma}

\def\E{\textbf{E}}

\begin{document}
\pagenumbering{arabic}
  \thispagestyle{empty}
   \chapter{On a Distribution Function Arising in Computational Biology}
\chapterauthors{Craig A.~Tracy and Harold Widom}

{\renewcommand{\thefootnote}{\fnsymbol{footnote}}

\footnotetext{\kern-19pt{\bf AMS Subject classification:} Primary
05A16; secondary 60F10, 92D20. }

\footnotetext{\kern-19pt{\bf  Keywords and phrases:} 
computational biology, random permutations,
increasing subsequences, Painlev\'e II.
}

\begin{abstract}
Karlin and Altschul in their statistical analysis for multiple high-scoring segments in molecular sequences   
 introduced a distribution function which gives the probability there are at least $r$ distinct
 and consistently ordered segment
     pairs all with score at least $x$. For 
long sequences this distribution can be expressed in terms of the distribution of the length of
     the longest increasing subsequence in a random permutation. Within the past few years, 
this last quantity has been
     extensively studied in the mathematics literature. 
The purpose of this note is to summarize these new mathematical
     developments in a form suitable for use in computational biology.
\end{abstract}
\begin{center}
\textit{Dedicated to Barry McCoy on the occasion of his sixtieth birthday.}
\end{center}
\section{The Distribution  Function}\label{s:intro-fbm}
Karlin and Altschul \cite{karlin} in their statistical analysis
for multiple high-scoring segments in molecular sequences,
introduced the following distribution function:  Let $F(r;y)$
denote the probability that there are at least $r$ distinct
\textit{and consistently ordered} segment pairs all with score
at least $x$. They further introduced a parameter
$y=KNe^{-\la x}$ where $K$ and $\la$ are parameters
related to the scoring system,  see \cite{karlin} for details.  We use the parameter
$y$ without further reference to $x$.  For long sequences ($N\ra\iy$)
this distribution function is well approximated by \cite{karlin}
\begin{equation}
F(r;y)=e^{-y}\,\sum_{k=r}^{\iy}{y^k R_{k,r}\ov k!^2},
\> r=1,2,\ldots,  \label{distrFn1} \end{equation}
where $R_{k,r}$ is the number of permutations of the integers $\{1,\ldots,k\}$
that contain an increasing subsequence of length at least $r$.
Let $X_y$ denote a positive integer valued random variable such that
\[ \pr\left(X_y\ge r\right)=F(r;y). \]

If $R^c_{k,r}$ denotes the complement of $R_{k,r}$, i.e.\ the 
number of
permutations $\s\in S_k$ all of whose increasing subsequences have
length strictly less than $r$, then clearly
\ba R_{k,r}^c &=& \#\left\{\s\in S_k: \ell_k(\s)<r \right\} \\
&=&\#\left\{\s\in S_k: \ell_k(\s)\le r-1 \right\}\\
&:=& f_{k,r-1}\ea
where $\ell_k(\s)$ is the length of the longest increasing
subsequence in $\s\in S_k$. 

\Remarks
\vspace{-1ex}
\begin{enumerate}
\item $F(r;y)$ is a distribution function in $r$ with \textit{parameter} $y$.
\item Dropping the requirement of consistent ordering has the effect
of replacing $R_{k,r}$ by
$k!$ in (\ref{distrFn1}).
 Thus the segments are Poisson distributed with parameter $y$.
\end{enumerate}

\section{Summary of Known Properties}
By convention, $f_{0,r}:=1$ for all $r$ and we note that
$f_{k,r}=k!$ if $k\le r$.  It is also convenient to introduce
the parameter $t\ge 0$, 
\[ y= t^2. \]
If we define
\begin{equation}
 D_r(t)=\sum_{k=0}^{\iy} {f_{k,r} t^{2k}\ov k!^2}\,, \label{distrFn2}\end{equation}
then
\bae
F(r;y)&=&e^{-y}\sum_{k=r}^{\iy}\left\{ k! - f_{k,r-1}\right\} \nonumber \\
      &=& e^{-y}\sum_{k=r}^{\iy} {y^k\ov k!} - e^{-y} \sum_{k=r} {
      f_{k,r-1} y^k\ov k!^2} \nonumber \\
      &=&e^{-y}\sum_{k=r}^{\iy} {y^k\ov k!}-e^{-y}\left(D_{r-1}(\sqrt{y})
      -\sum_{k=0}^r{y^k\ov k!}\right)\nonumber \\
      &=& 1- e^{-y} D_{r-1}(\sqrt{y}).\label{distrFn3}
      \eae
      
From Gessel~\cite{gessel} we know that $D_r(t)$ is the $r\times r$
Toeplitz determinant with symbol
\[ f(z)= e^{t(z+1/z)}. \]
In the past few years, $D_r$ has been extensively studied in
connection with the limiting distribution of the length of the longest
increasing subsequence of a random permutation, see Baik, Deift
and Johansson \cite{bdj} and
Aldous and Diaconis \cite{aldous1, aldous2}. 
We now summarize some of the these results.  Gessel's theorem
says that for all $r=1,2,\ldots$
\[ D_r(t)=\det\left(b_{\vert i-j\vert}\right)_{0\le i,j \le r-1} \]
where
$b_j:=I_j(2t)$
and $I_j$ is the modified Bessel function.  For \textit{small}
$r$ one simply evaluates this determinant to obtain
\ba
F(1;y)&=&1-e^{-y}, \\
F(2;y)&=&1-b_0 \,e^{-y}, \\
F(3;y)&=&1-\left(b_0^2-b_1^2\right)\, e^{-y},\\
F(4;y)&=&1-\left(b_0^3+2 b_1^2 b_2-2 b_1^2 b_0-b_0 b_2^2\right)\, e^{-y}.
\ea
From (\ref{distrFn3}) we see that
\[ \phi_r(y):= e^{-y} D_r(\sqrt{y})=\pr\left(X_y\le r\right). \]
Johansson~\cite{johansson} has shown that
for any given $\ve>0$, there exist $C$ and $\dl>0$ such that
\ba
0\le \phi_r(y) \le C e^{-\dl y} & \quad\textrm{if}\quad & (1+\ve)r<2\sqrt{y}, \\
0\le 1-\phi_r(y) \le {C\ov r} & \quad\textrm{if}\quad & (1-\ve)r>2\sqrt{y}. 
\ea
The breakthrough result of Baik-Deift-Johansson \cite{bdj} is the sharper asymptotic result
\begin{equation}
\lim_{y\ra\iy} \phi_{2\sqrt{y} + s \,y^{1/6}}(y)=F_2(s) \label{bdj}
\end{equation}
where $F_2$ is the  distribution function, first discovered by the present
authors in the context of random matrix theory \cite{tw0, tw1}
(see \cite{tw3} for a review), 
\begin{equation}
F_2(s)=\exp\left(-\int_s^{\iy} (x-s) q(x)^2\,dx\right)\label{F2}\end{equation}
and $q$ is the solution of the Painlev\'e II equation
\begin{equation}
q''=s q + 2 q^3 \label{p2} \end{equation}
satisfying $q(s)\sim \textrm{Ai}(s)$ as $s\ra\iy$.  (Here $\textrm{Ai}$
is the Airy function.)  It is known that such a solution
to (\ref{p2}) exists and is unique.
 A graph
of the density $dF_2/ds$  as well as some
 statistics of $F_2$ can be found in \cite{tw4}.
In terms of the random variable $X_y$ this says

\[ \chi_y:= {X_y-2\sqrt{y}\ov y^{1/6}} \]
converges weakly to a random variable, call it $\chi$, with distribution function
$F_2$. It was also proved that the scaled moments converge to the moments
of $F_2$ \cite{bdj}.

For finite $r$ we now describe some results of 
Periwal and Shevitz \cite{periwal}, Hisakado \cite{hisakado},
 Tracy and Widom \cite{tw2},
and Adler and van Moerbeke \cite{adlervanM}.  (We follow
the notation of \cite{tw2}.)  We have the representation
\begin{equation}
\phi_r(y) = \exp\left(-4\int_0^t \log(t/\t)\, \t\, (1-\p(\t))\,d\t\right),
\>\> y=t^2,\label{P5Rep}
\end{equation}
where $\p$ as a function of $t$ satisfies the equation
\begin{equation}
\p''={1\ov 2}\left({1\ov \p-1}+{1\ov\p}\right)(\p')^2-{1\ov t}\p'-8\p(\p-1)
+ 2\, {r^2\ov t^2}\,{\p-1\ov \p}\, . \label{phiDE}\end{equation}
We want the solution $\p$ that satisfies
\begin{equation} \p=1- {t^{2r}\ov (r!)^2}+\textrm{O}(t^{2r+1}), \>\> t\ra 0. \label{phiBC}
\end{equation}
Setting
\[ U_r^2:=1-\p, \]
we have the recursion relation, sometimes referred to as the discrete Painlev\'e II
equation,
\begin{equation} {r\ov t}\, U_r + (1-U_r^2)(U_{r-1}+U_{r+1}) = 0,
\>\> r=1,2,\ldots. \label{Urecursion}\end{equation}
The initial conditions for this recursion relation can be obtained from
$\phi_0=e^{-y}$  and $\phi_1=b_0 e^{-y}$.  A computation shows\footnotemark[2]
\footnotetext[2]{The
signs of $U_0$ and $U_1$ are not fixed from $\phi_0$ and $\phi_1$.
In \cite{tw2} the leading small $t$ behavior of $U_n$ is computed.  We
use this to fix the signs of $U_0$ and $U_1$.}
\[ U_0=-1, \>\> U_1={I_1(2t)\ov I_0(2t)}. \]

To make computational use of 
this distribution function, one needs  \textit{computationally} feasible
formulas for the first few moments of $X_y$ for all $y$; and more generally,
the distribution function itself.  Here are some partial results.
Of course,
\begin{equation}
 \E(X_y)=\sum_{r=1}^{\iy}r\left(\phi_r-\phi_{r-1}\right).
 \label{firstMoment}\end{equation} 
From (\ref{P5Rep}) and (\ref{phiBC}) it follows that
\[ \phi_r(y)=1-{y^{r+1}\ov (r+1)!^2} +\textrm{O}(y^{r+2}) \]
and thus
\[ \E(X_y)=\sum_{r=1}^Rr\left(\phi_r-\phi_{r-1}\right) +\textrm{O}(y^{R+1}). \]
Using the recursion relation (\ref{Urecursion}) we can  compute
the first few $U_n$'s and expand these for small $y$.  In this 
way we derive
\begin{equation}
 \E(X_y)=y -{1\ov 4}\, y^2 +{1\ov 12}\, y^3 - {7\ov 288}\, y^4 + \textrm{O}(y^5)
\label{mu1Small}\end{equation}
and similarly for the variance
\begin{equation}
 \textrm{Var}(X_y)=y-{3\ov 4}\,y^2+{17\ov 36}\, y^3 - {67\ov 288}\, y^4
+ \textrm{O}(y^5). \label{varSmall}\end{equation}
Note that the leading order terms are Poisson.
Higher order expansion coefficients  are given in Table~\ref{smallYdata}.

From \cite{bdj} we know that for large $y$, the
essential contribution in (\ref{firstMoment}) comes from $r$
around $2\sqrt{y}$.  Thus  for  $y\ra\iy$
\bae \E(X_y)&=&2\sqrt{y}+y^{1/6}\, \E(\chi)+\textrm{o}(y^{1/6}) \label{muLarge1} \\
&= & 2\sqrt{y} -1.77109\, y^{1/6} +\textrm{o}(y^{1/6}) \label{mu1Large2}\eae
where $\chi$ has distribution function (\ref{F2}).
We note for future reference, 
\[ \textrm{Var}(\chi)\approx 0.8132. \]
The small $y$ expansion of $\E(X_y)$ was computed through order 20.  If
we demand that the last coefficient in this expansion be less than, say, $1/10$, then $y<7.8$. 
 Evaluating this expansion at $y=7.8$ gives $\E(X_y)=3.66$ whereas the large $y$ expansion evaluated at
$y=7.8$ equals $3.09$ which is a difference of $0.57$.
To improve the overlap of these two expansions, one needs to compute the error
term in (\ref{muLarge1}).  

\section{An Example}
Karlin and Altschul give the parameters in their theory for the pairwise sequence comparison
of the chicken gene X protein and the fowlpox virus antithrombin III homolog.  The scoring
system gives $\lambda=0.314$, $K=0.17$ and $N=34336$.  For the three alignments found 
(see Table 3 in \cite{karlin})
the normalized scores (values of $x$) are 7.6, 6.7 and 5.8.  Using the above distribution
function we compute the expected number of distinct consistently ordered
seqment pairs with at least normalized score $x$.
The results are displayed in Table~\ref{bioData}.

\begin{table}
\begin{center}
{\large
\begin{tabular}{|c|c|c|c|}\hline 
$x$ & $y$ & $\E(X_y)$ & $\textrm{Var}(X_y)$ \\[0.5ex]\hline
7.6 & 536.8 & 41.3 & 6.6 \\[0.5ex]
6.7 & 712.1 & 48.1 & 7.3 \\[0.5ex]
5.8 & 944.6 & 55.9 & 8.0 \\[0.5ex]
\hline
\end{tabular}
}
\vspace{2ex}
\caption{\label{bioData}
For Karlin-Altschul parameters $\lambda=0.314$,
$K=0.17$ and $N=34336$, and three values of the
normalized score $x$, the expected value and variance
of $X_y$ are computed using the large $y$ expansions.
}
\end{center}
\end{table}

\begin{table}
\begin{center}
{\Large
\begin{tabular}{|r|r|r|}\hline 
$r$ & $c_r^1$ & $c_r^2$\\[0.5ex]
\hline
1 & 1 & 1 \\
2 & $-{1\ov 2^2}$ & $-{3\ov 2^2}$ \\[0.5ex]
3 & ${1\ov 2^2\, 3}$ & ${17\ov 2^2\, 3^2}$\\[0.5ex]
4 & $-{7\ov 2^5\, 3^2}$ & $-{67\ov 2^5\, 3^2}$ \\ [0.5ex]
5 & ${17\ov 2^6 \,3^2\, 5}$ & ${269\ov 2^6\,3^2\, 5}$ \\[0.5ex]
6 & $-{619\ov 2^8 \,3^4 \,5^2 }$ & $-{13\,\cdot\,
19\,\cdot\, 67 \ov 2^8\, 3^4\, 5^2}$\\[0.5ex]
7 & ${41\ov 2^7\, 3^2 \,5^2\, 7}$ & ${3491\ov 2^7\, 3^4\, 5^1\, 7}$ \\[0.5ex]
8 & $-{4001\ov 2^{12}\, 3^3\, 5^2\, 7^2}$ &$-{1064243\ov 2^{12}\, 3^4\, 5^2\, 7^2}$\\[0.5ex]
9 & ${173\,\cdot\,313\ov 2^{14}\, 3^6\, 5^2\, 7^2}$ &
${28638487\ov 2^{14}\, 3^7\, 5^2\, 7^2}$\\[0.5ex]
10 & $-{17\,\cdot\,62687\ov 2^{16}\,3^8\,5^3\,7^2}$ & 
$-{41\,\cdot\,557\,\cdot\,17257\ov 2^{16}\, 3^8\, 5^3\, 7^2}$ \\[0.5ex]
11& ${2823631\ov 2^{15}\,3^8\,5^4\,7^2\,11}$ &
${37\,\cdot\, 61924123\ov 2^{15}\,3^8\, 5^4\, 7^2\, 11}$\\[0.5ex]
12 & $-{941\,\cdot\, 407219\ov 2^{19}\, 3^{10}\, 5^4\, 7^2\, 11^2}$ &
$-{17\,\cdot\,29\,\cdot 286954607\ov 2^{19}\, 3^{10}\, 5^3 \, 7^2\, 11^2}$\\[0.5ex]
13 & ${6377893 \ov 2^{17}\, 3^{9}\, 5^3\, 7^2\, 11^2 \, 13}$ &
${206619709873\ov 2^{16}\, 3^{10} \, 5^4\, 7^2\, 11^2\, 13}$\\[0.5ex]
14 & $-{11657\,\cdot\, 1658989\ov 2^{22}\,3^{10}\,5^4\,7^3\,11^2\,13^2}$
& $-{199735173503123\ov 2^{22}\, 3^{10}\, 5^4\, 7^3\, 11^2\, 13^2}$  \\[0.5ex]
15 & ${179\,\cdot\,257\,\cdot\,139493\ov 2^{20}\, 3^{11}\, 5^4\, 7^4\, 11^2\,
13^2}$ & ${479147\,\cdot\, 50402324263\ov 2^{21}\, 3^{12}\, 5^6\, 7^4\, 11^2\,
13^2}$\\[0.5ex]
16& $-{37\,\cdot\,23593\,\cdot\,1363963\ov 2^{27}\,3^{11}\,5^6\,7^4\,11^2\,13^2}$ &
$-{59\,\cdot\,163363\,\cdot\,7608612619\ov 2^{27}\, 3^{11}\, 5^6\, 7^4\, 11^2\, 13^2}$\\[0.5ex]
17& ${43\,\cdot\,863\,\cdot\,701781161\ov 2^{30}\, 3^{12}\, 5^6\, 7^4\, 11^2\, 13^2\, 17}$ &
${27057479\,\cdot\,146285342603\ov 2^{30}\, 3^{12}\, 5^6\, 7^4\, 11^2\, 13^2\, 17}$\\[0.5ex]
18 & $-{23\,\cdot\,5264671\,\cdot\,6578291\ov 2^{32}\, 3^{14}\, 5^6\, 7^4\, 11^2\, 13^2\, 17^2}$&
$-{307\,\cdot\,972530242052278499\ov 2^{32}\, 3^{14}\, 5^6\, 7^4\, 11^2\, 13^2\, 17^2}$\\[0.5ex]
19&${1077161\,\cdot\,39636029\ov 2^{31}\, 3^{15}\, 5^4\, 7^4\, 11^2\, 13^2\, 17^2\, 19}$&
${61\,\cdot\,83\,\cdot\,709\,\cdot\,7309\,\cdot\,37338914351\ov
2^{31}\,3^{15}\,5^6\,7^4\,11^2\,13^2\,17^2\,19}$\\[0.5ex]
20&$-{229\,\cdot\,5189\,\cdot\,247913\,\cdot\,1229957\ov
2^{35}\, 3^{15}\, 5^8\, 7^4\, 11^2\, 13^2\, 17^2\, 19^2}$ &
$-{239\,\cdot\,1181\,\cdot\,2161\,\cdot\,263188412702251\ov
2^{35}\, 3^{15}\, 5^7\, 7^4\, 11^2\, 13^2\, 17^2\, 19^2}$ \\[0.5ex]
\hline
\end{tabular}
}
\vspace{2ex}
\caption{\label{smallYdata}
The number $c_r^1$  is the coefficient of $y^r$ in the
small $y$ expansion of $\E(X_y)$ and $c_r^2$ is
the coefficient of $y^r$ in the small $y$
expansion of $\textrm{Var}(X_y)$.}
\end{center}
\end{table}

\Ack The authors thank Professor S.~Karlin and Dr.~S.~Altschul for helpful
comments and for allowing us to
use their data in Table \ref{bioData}.
The first author wishes to thank Professors M.~Kashiwara and T.~Miwa
for the invitation to speak at \textit{MathPhys Odyssey 2001: Integrable
Models and Beyond}. 
This work was partially supported by the NSF Grants
 DMS--9802122 and DMS--9732687.

\newpage


\vskip 2pc
{\noindent Craig A.~Tracy,
 Department of Mathematics and
 Institute of Theoretical Dynamics,
 University of California,
 Davis, CA 95616,    U.S.A.,
 email: {\tt tracy@itd.ucdavis.edu }
 \vskip2ex
\noindent Harold Widom,
Department of Mathematics,
University of California,
Santa Cruz, CA 95064,
email: {\tt widom@math.ucsc.edu }}
\end{document}